\documentclass[11pt]{amsart}
\usepackage{amssymb}
\usepackage{amsmath}
\usepackage[active]{srcltx}
\usepackage{t1enc}
\usepackage[latin2]{inputenc}
\usepackage{verbatim}
\usepackage{amsmath,amsfonts,amssymb,amsthm}
\usepackage[mathcal]{eucal}
\usepackage{enumerate}
\usepackage[centertags]{amsmath}
\usepackage{graphics}

\newtheorem{theorem}{Theorem}

\newtheorem{corollary}{Corollary}

\newtheorem{proposition}{Proposition}
\begin{document}
\author{Giorgi Tutberidze}
\title[ $T$ means]{Maximal operators of $T$ means with respect to the Vilenkin system }
\address{G.Tutberidze, The University of Georgia, School of science and technology, 77a Merab Kostava St, Tbilisi 0128, Georgia and Department of Computer Science and Computational Engineering, UiT - The Arctic University of Norway, P.O. Box 385, N-8505, Narvik, Norway.}
\email{giorgi.tutberidze1991@gmail.com}
\thanks{The research was supported by Shota Rustaveli National Science Foundation grant FR-19-676.}
\date{}
\maketitle

\begin{abstract}
In this paper we prove and discuss some new $\left( H_{p},weak-L_{p}\right) $
type inequalities of maximal operators of $T$ means with respect to Vilenkin systems with monotone coefficients. We also apply these results to prove a.e. convergence
of \ such $T$ means. It is also proved that these results are the best possible in a special sense. As applications, both some well-known and new results are pointed out.
\end{abstract}



\section{Introduction}

For the notation used in this introduction see Section 2.

\bigskip\ Weisz
\cite{We2} proved boundedness of $\ \sigma
^{\ast }$ \ from the martingale space $H_{p}$ to the space $L_{p},$ for $%
p>1/2$. Simon \cite{Si1} gave a counterexample, which shows that boundedness
does not hold for $0<p<1/2.$ A counterexample for $p=1/2$ was given by
Goginava \cite{Go} (see also \cite{tep1, tep2} and \cite{pttw}). Moreover, Weisz \cite{We4} proved
that the following is true:

\textbf{Theorem W1. }The maximal operator of the Fejér means $\sigma ^{\ast
} $ is bounded from the Hardy space $H_{1/2}$ to the space weak-$L_{1/2}$.

Riesz`s logarithmic means with respect to the Walsh and
Vilenkin systems were investigated by Simon \cite{Si1},
Blahota and Gát \cite{bg}. For the Vilenkin systems in \cite{tep5} and for the Walsh system in \cite{sw} it were proved that the maximal operator of Riesz`s means $R^{\ast }$ is bounded from the Hardy space $H_{1/2}$ to the space
$weak-L_{1/2}$, but is not bounded from the Hardy space $H_{p}$ to the space $L_{p},$ when $0<p\leq 1/2.$ Since the set of Vilenkin polynomials are dense in $L_{1},$ by well-known density argument due to Marcinkiewicz and Zygmund \cite{mz} we have that $R_{n}f\rightarrow f,$ \ a.e. for \ all \ $f\in
L_{1} $.

Móricz and Siddiqi \cite{Mor} investigated the approximation properties of some special Nörlund means of Walsh-Fourier series of $L_{p}$ function in norm. In the two-dimensional case similar problems was studied by Nagy  \cite{n,nagy}. In \cite{ptw} (see also \cite{BPTW, bpt1}) it was proved some $(H_{p},L_{p})$-type inequalities for the maximal operators of Nörlund means, when $0<p\leq 1.$

In \cite{BN} and \cite{bnt1} were investigated $T$ means and studied some approximation properties of these summability methods in the Lebesgue spaces for $L_p,$ $1\leq p \leq \infty.$ In this paper we prove analogous theorems considered in \cite{ptw} and derive some new $(H_{p},L_{p})$-type inequalities for the maximal operators of $T$ means, when $0<p\leq 1.$ We also apply these results to prove a.e. convergence of \ such $T$ means. It is also proved that these results are the best possible in a special sense. As applications, both some well-known and new results are pointed out.

The paper is organized as follows: In Section 3 we present and discuss the main results and in Section 4 the proofs can be found. Moreover, in order not to disturb our discussions in these Sections some preliminaries are given in Section 2.

\section{Preliminaries}

Denote by $\mathbb{N}_{+}$ the set of the positive integers, $\mathbb{N}:=\mathbb{N}_{+}\cup \{0\}.$ Let $m:=(m_{0,}$ $m_{1},...)$ be a sequence of the positive
integers not less than 2. Denote by
\begin{equation*}
Z_{m_{k}}:=\{0,1,...,m_{k}-1\}
\end{equation*}
the additive group of integers modulo $m_{k}$.

Define the group $G_{m}$ as the complete direct product of the groups $%
Z_{m_{i}}$ with the product of the discrete topologies of $Z_{m_{j}}`$s.

The direct product $\mu $ of the measures
\begin{equation*}
\mu _{k}\left( \{j\}\right) :=1/m_{k}\text{ \ \ \ }(j\in Z_{m_{k}})
\end{equation*}%
is the Haar measure on $G_{m_{\text{ }}}$with $\mu \left( G_{m}\right) =1.$

In this paper we discuss bounded Vilenkin groups,\textbf{\ }i.e. the case
when $\sup_{n}m_{n}<\infty .$

The elements of $G_{m}$ are represented by sequences
\begin{equation*}
x:=\left( x_{0},x_{1},...,x_{j},...\right) ,\ \left( x_{j}\in
Z_{m_{j}}\right) .
\end{equation*}

It is easy to give a base for the neighborhood of $G_{m}:$

\begin{equation*}
I_{0}\left( x\right) :=G_{m},\text{ \ }I_{n}(x):=\{y\in G_{m}\mid
y_{0}=x_{0},...,y_{n-1}=x_{n-1}\},
\end{equation*}%
where $x\in G_{m},$ $n\in
\mathbb{N}
.$

Denote $I_{n}:=I_{n}\left( 0\right) $ for $n\in
\mathbb{N}
_{+},$ and $\overline{I_{n}}:=G_{m}$ $\backslash $ $I_{n}$.

\bigskip If we define the so-called generalized number system based on $m$
in the following way :
\begin{equation*}
M_{0}:=1,\ M_{k+1}:=m_{k}M_{k}\,\,\,\ \ (k\in
\mathbb{N}
),
\end{equation*}%
then every $n\in
\mathbb{N}
$ can be uniquely expressed as $n=\sum_{j=0}^{\infty }n_{j}M_{j},$ where $%
n_{j}\in Z_{m_{j}}$ $(j\in
\mathbb{N}
_{+})$ and only a finite number of $n_{j}`$s differ from zero.

Next, we introduce on $G_{m}$ an orthonormal system which is called the
Vilenkin system. At first, we define the complex-valued function $%
r_{k}\left( x\right) :G_{m}\rightarrow
\mathbb{C}
,$ the generalized Rademacher functions, by%
\begin{equation*}
r_{k}\left( x\right) :=\exp \left( 2\pi ix_{k}/m_{k}\right) ,\text{ }\left(
i^{2}=-1,x\in G_{m},\text{ }k\in
\mathbb{N}
\right) .
\end{equation*}

Now, define the Vilenkin system$\,\,\,\psi :=(\psi _{n}:n\in
\mathbb{N}
)$ on $G_{m}$ as:
\begin{equation*}
\psi _{n}(x):=\prod\limits_{k=0}^{\infty }r_{k}^{n_{k}}\left( x\right)
,\,\,\ \ \,\left( n\in
\mathbb{N}
\right) .
\end{equation*}

Specifically, we call this system the Walsh-Paley system when $m\equiv 2.$

The norms (or quasi-norms) of the spaces $L_{p}(G_{m})$ and $weak-L_{p}\left( G_{m}\right) $ $\left( 0<p<\infty
\right) $ are respectively defined by
\begin{equation*}
\left\Vert f\right\Vert _{p}^{p}:=\int_{G_{m}}\left\vert f\right\vert ^{p}d\mu \ \ \ \text{and} \ \ \ \left\Vert f\right\Vert _{weak-L_{p}}^{p}:=\underset{\lambda >0}{\sup }
\lambda ^{p}\mu \left( f>\lambda \right) <+\infty .
\end{equation*}%

The Vilenkin system is orthonormal and complete in $L_{2}\left( G_{m}\right)
$ (see \cite{Vi}).

Now, we introduce analogues of the usual definitions in Fourier-analysis. If
$f\in L_{1}\left( G_{m}\right) $ we can define the Fourier coefficients, the
partial sums of the Fourier series, the Dirichlet kernels with respect to
the Vilenkin system in the usual manner:

\begin{equation*}
\widehat{f}\left( n\right) :=\int_{G_{m}}f\overline{\psi }_{n}d\mu \,\ \ \ \
\ \ \,\left( n\in
\mathbb{N}
\right) ,
\end{equation*}%
\begin{equation*}
S_{n}f:=\sum_{k=0}^{n-1}\widehat{f}\left( k\right) \psi _{k},\text{ \ \ }%
D_{n}:=\sum_{k=0}^{n-1}\psi _{k\text{ }},\text{ \ \ }\left( n\in
\mathbb{N}
_{+}\right)
\end{equation*}%
respectively. Recall that
\begin{equation}
D_{M_{n}}\left( x\right) =\left\{
\begin{array}{ll}
M_{n}, & \text{if\thinspace \thinspace \thinspace }x\in I_{n}, \\
0, & \text{if}\,\,x\notin I_{n}.%
\end{array}%
\right.  \label{1dn}
\end{equation}

The $\sigma $-algebra generated by the intervals $\left\{ I_{n}\left(
x\right) :x\in G_{m}\right\} $ will be denoted by $\digamma _{n}\left( n\in\mathbb{N} \right) .$ Denote by $f=\left( f^{\left( n\right) },n\in\mathbb{N} \right) $ a martingale with respect to $\digamma _{n}\left( n\in\mathbb{N}\right)$ (for details see e.g. \cite{We1}).

For $0<p<\infty $ \ the Hardy martingale spaces $H_{p}\left( G_{m}\right) $
consist of all martingales for which
\begin{equation*}
\left\Vert f\right\Vert _{H_{p}}:=\left\Vert f^{\ast }\right\Vert
_{p}<\infty, \ \ \ \text{where}\ \ \ f^{\ast }:=\sup_{n\in\mathbb{N}}\left\vert f^{(n)}\right\vert .
\end{equation*}

A bounded measurable function $a$ is called a p-atom, if there exists an interval $I$, such that
\begin{equation*}
\int_{I}ad\mu =0,\text{ \ \ }\left\Vert a\right\Vert _{\infty }\leq \mu
\left( I\right) ^{-1/p},\text{ \ \ supp}\left( a\right) \subset I.
\end{equation*}

Weisz \cite{We3} proved that the Hardy spaces $H_p$ have atomic characterizations. In particular the following is true:
\begin{proposition}\label{prop1} A martingale $f=\left( f^{\left( n\right) },n\in	\mathbb{N}	\right) $ is in $H_{p}\left( 0<p\leq 1\right) $ if and only if there exists a sequence $\left( a_{k},k\in	\mathbb{N} \right) $ of p-atoms and a sequence $\left( \mu _{k},k\in	\mathbb{N} \right) ,$ of real numbers, such that, for every $n\in\mathbb{N},$
\end{proposition}
\begin{equation}
\qquad \sum_{k=0}^{\infty }\mu _{k}S_{M_{n}}a_{k}=f^{\left( n\right) },\ \ \ \ \ \ \ \sum_{k=0}^{\infty }\left\vert \mu _{k}\right\vert ^{p}<\infty .
\label{1}
\end{equation}
Moreover,
\begin{equation*}
\left\Vert f\right\Vert _{H_{p}}\backsim \inf \left( \sum_{k=0}^{\infty
}\left\vert \mu _{k}\right\vert ^{p}\right) ^{1/p},
\end{equation*}
\textit{where the infimum is taken over all decomposition of} $f$ \textit{of the form} (\ref{1}).
We also need the following result of Weisz \cite{We3}:
\begin{proposition}\label{prop2} Suppose that the operator $T$ \ is $\sigma $-linear and for some $
	0<p<1$ 
	\begin{equation*}
	\left\Vert Tf\right\Vert _{weak-L_{p}}\leq c_{p}\left\Vert f\right\Vert
	_{H_{p}},
	\end{equation*}
	then $T$ is of weak	type-(1,1) i.e.
	\begin{equation*}
	\left\Vert Tf\right\Vert _{weak-L_{1}}\leq c\left\Vert f\right\Vert _{1}.
	\end{equation*}
\end{proposition}

If $f=\left( f^{\left( n\right) },n\in\mathbb{N}\right) $ is a martingale, then the Vilenkin-Fourier coefficients must be
defined in a slightly different manner:
\begin{equation*}
\widehat{f}\left( i\right) :=\lim_{k\rightarrow \infty
}\int_{G_{m}}f^{\left( k\right) }\overline{\psi }_{i}d\mu .
\end{equation*}

Let $\{q_{k}:k\geq 0\}$ be a sequence of non-negative numbers. The $n$-th Nörlund and  $T$ means for a Fourier series of $f$ \ are respectively defined by

\begin{equation*}
t_{n}f=\frac{1}{Q_{n}}\overset{n}{\underset{k=1}{\sum }}q_{n-k}S_{k}f,
\end{equation*}
and
\begin{equation} \label{nor}
T_nf:=\frac{1}{Q_n}\overset{n-1}{\underset{k=0}{\sum }}q_{k}S_kf,
\end{equation}
where $Q_{n}:=\sum_{k=0}^{n-1}q_{k}.$ 
It is obvious that 
\begin{equation*}
	T_nf\left(x\right)=\underset{G_m}{\int}f\left(t\right)F_n\left(x-t\right) d\mu\left(t\right),
\end{equation*}
where $	F_n:=\frac{1}{Q_n}\overset{n}{\underset{k=1}{\sum }}q_{k}D_k$ is called the $T$ kernel.

We always assume that $\{q_k:k\geq 0\}$ is a sequence of non-negative numbers and $q_0>0.$ Then the summability method (\ref{nor}) generated by $\{q_k:k\geq 0\}$ is regular if and only if $	\lim_{n\rightarrow\infty}Q_n=\infty.$

If we invoke Abel transformation we get the following identities, which are very important for the investigations of $T$ summability:
\begin{eqnarray} \label{2b}
Q_n&:=&\overset{n-1}{\underset{j=0}{\sum}}q_j =\overset{n-2}{\underset{j=0}{\sum}}\left(q_{j}-q_{j+1}\right) j+q_{n-1}{(n-1)}
\end{eqnarray}
and
\begin{equation} 	\label{2c}
F_n=\frac{1}{Q_n}\left(\overset{n-2}{\underset{j=0}{\sum}}\left(q_j-q_{j+1}\right) jK_j+q_{n-1}(n-1)K_{n-1}\right).
\end{equation}

The well-known example of N\"orlund summability is the so-called $\left(C,\alpha\right)$-mean (Ces\`aro means), which are defined by
\begin{equation*}
	\sigma_n^{\alpha}f:=\frac{1}{A_n^{\alpha}}\overset{n}{\underset{k=1}{
			\sum}}A_{n-k}^{\alpha-1}S_kf, \qquad 0<\alpha<1,
\end{equation*}
where 
\begin{equation*}
	A_0^{\alpha}:=0,\qquad A_n^{\alpha}:=\frac{\left(\alpha+1\right)...\left(\alpha+n\right)}{n!},\qquad \alpha \neq -1,-2,...
\end{equation*}

We also consider the "inverse" $\left(C,\alpha\right)$-means, which is an example of a $T$-means:
\begin{equation*}
	U_n^{\alpha}f:=\frac{1}{A_n^{\alpha}}\overset{n-1}{\underset{k=0}{\sum}}A_{k}^{\alpha-1}S_kf, \qquad 0<\alpha<1.
\end{equation*}

Let $V_n^{\alpha}$ denote
the $T$ mean, where $	\left\{q_0=1, \  q_k=k^{\alpha-1}:k\in \mathbb{N}_+\right\} ,$
that is 
\begin{equation*}
	V_n^{\alpha}f:=\frac{1}{Q_n}\overset{n}{\underset{k=1}{\sum }}k^{\alpha-1}S_kf,\qquad 0<\alpha<1.
\end{equation*}

The $n$-th Riesz`s logarithmic mean $R_{n}$ and the Nörlund logarithmic mean
$L_{n}$ are defined by
\begin{equation*}
	R_{n}f:=\frac{1}{l_{n}}\sum_{k=1}^{n-1}\frac{S_{k}f}{k}\text{ \ \ \ and  \ \ \ }
	L_{n}f:=\frac{1}{l_{n}}\sum_{k=1}^{n-1}\frac{S_{k}f}{n-k},
\end{equation*}%
\ respectively, where $l_{n}:=\sum_{k=1}^{n-1}1/k.$

Up to now we have considered $T$ means in the case when the sequence $\{q_k:k\in\mathbb{N}\}$ is bounded but now we consider $T$ summabilities with unbounded sequence $\{q_k:k\in\mathbb{N}\}$. Let $\alpha\in
\mathbb{R}_+,\ \ \beta\in\mathbb{N}_+$ and
\begin{equation*}
	\log^{(\beta)}x:=\overset{\beta\text{times}}{\overbrace{\log ...\log}}x.
\end{equation*}

If we define the sequence $\{q_k:k\in \mathbb{N}\}$ by
$	\left\{q_0=0, \ q_k=\log^{\left(\beta \right)}k^{\alpha
	}:k\in\mathbb{N}_+\right\},$ 
then we get the class of $T$ means with non-decreasing coefficients:
\begin{equation*}
	B_n^{\alpha,\beta}f:=\frac{1}{Q_n}
	\sum_{k=1}^{n}\log^{\left(\beta\right)}k^{\alpha}S_kf.
\end{equation*}%

We note that $B_n^{\alpha,\beta}$ are
well-defined for every $n \in \mathbb{N}$

\begin{equation*}
	B_n^{\alpha,\beta}f=\sum_{k=1}^{n}\frac{\log^{\left(\beta\right)}k^{\alpha }}{Q_n}S_kf.
\end{equation*}

It is obvious that $\frac{n}{2}\log^{\left(\beta \right)}\frac{n^{\alpha }}{2^{\alpha }}\leq Q_n\leq n\log^{\left(\beta\right)}n^{\alpha}.$ It follows that
\begin{eqnarray} \label{node00}
	\frac{q_{n-1}}{Q_n}\leq\frac{c\log^{\left(\beta\right)}n^{\alpha}}{n\log^{\left(\beta\right) }n^{\alpha}}= O\left(\frac{1}{n}\right)\rightarrow 0,\text{ \ as \ }n\rightarrow \infty.
\end{eqnarray}

We also define the maximal operator $T^*$ of $T$ means by 
\begin{eqnarray*}
	T^{\ast}f:=\sup_{n\in\mathbb{N}}\left\vert T_nf\right\vert.
\end{eqnarray*}

Some well-known examples of maximal operators of $T$ means are the maximal operator of Fejér $\sigma^*$ and Riesz $R^*$ logarithmic means, which are defined by: 
\begin{eqnarray*}
	\sigma^{\ast}f:=\sup_{n\in\mathbb{N}}\left\vert \sigma_{n}f\right\vert, \ \ \ \ \ R^{\ast}f:=\sup_{n\in\mathbb{N}}\left\vert R_{n}f\right\vert.
\end{eqnarray*}

We also define some new maximal operators $U^{\alpha,\ast},V^{\alpha,\ast },B^{\alpha,\beta,\ast }, \ \ (\alpha\in \mathbb{R_+}, \beta \in \mathbb{N_+})$ by:
\begin{eqnarray*}
	U^{\alpha,\ast}f:=\sup_{n\in
		\mathbb{N}}\left\vert U_n^{\alpha }f\right\vert,\ \ \  V^{\alpha,\ast }f:=\sup_{n\in\mathbb{N}}\left\vert V_{n}^{\alpha}f\right\vert, \ \ \ 
	 B^{\alpha,\beta,\ast }f:=\sup_{n\in\mathbb{N}}\left\vert B_{n}^{\alpha,\beta}f\right\vert.
\end{eqnarray*}

\section{The Main Results}

First we state our main result concerning the maximal operator of the
summation method (\ref{nor}), which we also show is in a sense sharp.

\begin{theorem}\label{theorem1}
a) The maximal operator $T^{\ast }$ of \ the summability method (\ref{nor})
with non-increasing sequence $\{q_{k}:k\geq 0\},$ is bounded from the Hardy
space $H_{1/2}$ to the space $weak-L_{1/2}.$

The statement in a) is sharp in the following sense:

\bigskip b) Let $0<p<1/2$ and\ $\{q_{k}:k\geq 0\}$ be a non-increasing sequence, satisfying the condition
\begin{equation}
\frac{q_{n+1}}{Q_{n+2}}\geq \frac{c}{n},\text{ \ \ }\left( c\geq 1\right) .
\label{cond1}
\end{equation}%
Then there exists a martingale $f\in H_{p},$ such that
\begin{equation*}
\underset{n\in	\mathbb{N}}{\sup }\left\Vert T_{n}f\right\Vert _{weak-L_{p}}=\infty .
\end{equation*}
\end{theorem}

A number of special cases of our results are of particular interest and give both well-known and new information. We just give the following examples of such $T$ means with non-increasing sequence $\{q_{k}:k\geq 0\}:$

\begin{corollary}\label{cor0}
	\textbf{\ }The maximal operators $U^{\alpha, *}$, $V^{\alpha, *}$ and $R^{*}$ are  bounded from the Hardy space $H_{1/2}$ to the space $%
	weak-L_{1/2}$ but are not bounded from $H_{p}$ to the space $weak-L_{p},$
	when $0<p<1/2.$
\end{corollary}

\begin{corollary}\label{cor2}
	Let $f\in L_{1}$ and $T_{n}$ be the $T$ means with non-increasing sequence $\{q_{k}:k\geq 0\}$. Then
$\	T_{n}f\rightarrow f,\text{ \ \ a.e., \ \ as \ }n\rightarrow \infty .$
\end{corollary}

\begin{corollary}\label{cor1}
	Let $f\in L_{1}$. Then%
	\begin{eqnarray*}
		R_{n}f &\rightarrow &f,\text{ \ \ \ a.e., \ \ \ \ as \ }n\rightarrow \infty, \\
		U_{n} ^{\alpha}f &\rightarrow &f,\text{ \ \ \ a.e., \ \ \ as \ }n\rightarrow
		\infty ,\text{\ \ \ }\\
		V_{n} ^ {\alpha}f &\rightarrow &f,\text{ \ \ \ a.e., \ \ \ as \ }n\rightarrow
		\infty ,\text{\ \ \ }
	\end{eqnarray*}
\end{corollary}

Our next main result reads:
\begin{theorem} \label{theorem2}
a) The maximal operator $T^{\ast }$ of \ the summability method (\ref{nor})
with non-decreasing sequence $\{q_{k}:k\geq 0\}$ satisfying the condition 
\begin{equation}
\frac{q_{n-1}}{Q_n}= O\left(\frac{1}{n}\right)\label{Tmeanscond}
\end{equation}
is bounded from the Hardy space $H_{1/2}$ to the space $weak-L_{1/2}.$

b) Let $0<p<1/2.$ For any non-decreasing sequence $\{q_{k}:k\geq 0\},$ there exists a martingale $f\in H_{p},$ such that
\begin{equation*}
\underset{n\in
	\mathbb{N}
}{\sup }\left\Vert T_{n}f\right\Vert _{weak-L_{p}}=\infty .
\end{equation*}

\end{theorem}

A number of special cases of our results are of particular interest and give both well-known and new information. We just give the following examples of such  $T$ means with non-decreasing sequence $\{q_{k}:k\geq 0\}:$

\begin{corollary} \label{cor4}
\textbf{\ }The maximal operator $B^{\alpha,\beta, *}$ is  bounded from the Hardy space $H_{1/2}$ to the space $
weak-L_{1/2}$ but is not bounded from $H_{p}$ to the space $weak-L_{p},$
when $0<p<1/2.$
\end{corollary}

\begin{corollary}\label{cor5}
Let $f\in L_{1}$ and $T_{n}$ be the $T$ means with \textit{non-decreasing sequence }$\{q_{k}:k\geq 0\}$ \textit{and} satisfying condition (\ref{Tmeanscond}). Then
\begin{equation*}
T_{n}f\rightarrow f,\text{ \ \ a.e., \ \ as \ }n\rightarrow \infty .
\end{equation*}
\end{corollary}

\begin{corollary}\label{cor6}
	Let $f\in L_{1}$. Then $ \ B_{n}^{\alpha,\beta}f \rightarrow f,\text{ \ \ \ a.e., \ \ \ as \ }n\rightarrow 	\infty.$
\end{corollary}

\section{Proofs}

\textbf{Proof of Theorem \ref{theorem1}} a). Let the sequence $\{q_{k}:k\geq 0\}$ be non-increasing. By combining (\ref{2b}) with (\ref{2c}) and using Abel transformation we get that

\begin{eqnarray*}
\left\vert T_{n}f\right\vert &\leq &\left\vert \frac{1}{Q_{n}}\overset{n-1}{%
\underset{j=1}{\sum }}q_{j}S_{j}f\right\vert \\
&\leq &\frac{1}{Q_{n}}\left( \overset{n-2}{\underset{j=1}{\sum }}\left\vert
q_{j}-q_{j+1}\right\vert j\left\vert \sigma _{j}f\right\vert
+q_{n-1}(n-1)\left\vert \sigma _{n}f\right\vert \right) \\
&\leq &\frac{1}{Q_{n}}\left( \overset{n-2}{\underset{j=1}{\sum }}\left(
q_{j}-q_{j+1}\right) j+q_{n-1}(n-1)\right) \sigma ^{\ast }f\leq \sigma ^{\ast}f
\end{eqnarray*}
so that
\begin{equation}
T^{\ast }f\leq \sigma ^{\ast }f.  \label{12aaaa}
\end{equation}

If we apply (\ref{12aaaa}) and Theorem W1 we can conclude that the maximal operators  $T^{\ast }$ of all $T$ means with non-increasing sequence $\{q_{k}:k\geq 0\},$ are bounded from the Hardy space $H_{1/2}$ to the space $weak-L_{1/2}.$ The proof of part a) of Theorem 1 is complete.

b) Let $0<p<1/2$ and $\left\{ \alpha _{k}:k\in
\mathbb{N}
\right\} $ be an increasing sequence of positive integers such that:\qquad
\begin{equation}
\sum_{k=0}^{\infty }1/\alpha _{k}^{p}<\infty ,  \label{2}
\end{equation}

\begin{equation}
\lambda \sum_{\eta =0}^{k-1}\frac{M_{\alpha _{\eta }}^{1/p}}{\alpha _{\eta }}%
<\frac{M_{\alpha _{k}}^{1/p}}{\alpha _{k}},  \label{3}
\end{equation}

\begin{equation}
\frac{32\lambda M_{\alpha _{k-1}}^{1/p}}{\alpha _{k-1}}<\frac{M_{\alpha
		_{k}}^{1/p-2}}{\alpha _{k}},  \label{4}
\end{equation}%
where $\lambda =\sup_{n}m_{n}.$

We note that such an increasing sequence $\left\{ \alpha _{k}:k\in
\mathbb{N}
\right\} $ which satisfies conditions (\ref{2})-(\ref{4}) can be constructed.

Let \qquad
\begin{equation}
f^{\left( A\right) }=\sum_{\left\{ k;\text{ }\lambda _{k}<A\right\} }\lambda
_{k}a_{k},  \label{55}
\end{equation}%
where

\begin{equation*}
\lambda _{k}=\frac{\lambda }{\alpha _{k}}  \ \ \ \  and \ \ \ \
a_{k}=\frac{M_{\alpha _{k}}^{1/p-1}}{\lambda }\left( D_{M_{\alpha
		_{k}+1}}-D_{M_{_{\alpha _{k}}}}\right) .  \label{77}
\end{equation*}

By using Proposition \ref{prop1}, it is easy to show that the martingale $\,f=\left( f^{\left( 1\right)
},f^{\left( 2\right) }...f^{\left( A\right) }...\right) \in H_{1/2}.$ Moreover, it is easy to show that

\begin{equation}
\widehat{f}(j)=\left\{
\begin{array}{ll}
\frac{M_{\alpha _{k}}^{1/p-1}}{\alpha _{k}},\, & \text{if \thinspace
	\thinspace }j\in \left\{ M_{\alpha _{k}},...,M_{\alpha _{k}+1}-1\right\} ,%
\text{ }k=0,1,2..., \\
0, & \text{if \ \thinspace \thinspace \thinspace }j\notin
\bigcup\limits_{k=1}^{\infty }\left\{ M_{\alpha _{k}},...,M_{\alpha
	_{k}+1}-1\right\} .%
\end{array}%
\right.  \label{6}
\end{equation}

We can write

\begin{equation}
T_{M_{\alpha _{k}}+2}f=\frac{1}{Q_{M_{\alpha _{k}}+2}}\sum_{j=0}^{M_{\alpha
		_{k}}}q_{j}S_{j}f+\frac{q_{M_{\alpha
			_k}+1}}{Q_{M_{\alpha _{k}}+2}}S_{M_{\alpha
		_{k}}+1}f:=I+II. \label{155aba}
\end{equation}

Let $M_{\alpha _{s}}\leq $ $j\leq M_{\alpha _{s}+1},$ where $s=0,...,k-1.$
Moreover,
\begin{equation*}
\left\vert D_{j}-D_{M_{_{\alpha _{s}}}}\right\vert \leq j-M_{_{\alpha
		_{s}}}\leq \lambda M_{_{\alpha _{s}}},\text{ \ \ }\left( s\in
\mathbb{N}
\right)
\end{equation*}%
so that, according to (\ref{1dn}) and (\ref{6}), we have that%
\begin{eqnarray}
&&\left\vert S_{j}f\right\vert =\left\vert \sum_{v=0}^{M_{\alpha _{s-1}+1}-1}%
\widehat{f}(v)\psi _{v}+\sum_{v=M_{\alpha _{s}}}^{j-1}\widehat{f}(v)\psi
_{v}\right\vert  \label{8} \\
&\leq &\left\vert \sum_{\eta =0}^{s-1}\sum_{v=M_{\alpha _{\eta
}}}^{M_{\alpha _{\eta }+1}-1}\frac{M_{\alpha _{\eta }}^{1/p-1}}{\alpha
	_{\eta }}\psi _{v}\right\vert +\frac{M_{\alpha _{s}}^{1/p-1}}{\alpha _{s}}%
\left\vert \left( D_{j}-D_{M_{_{\alpha _{s}}}}\right) \right\vert  \notag \\
&=&\left\vert \sum_{\eta =0}^{s-1}\frac{M_{\alpha _{\eta }}^{1/p-1}}{\alpha
	_{\eta }}\left( D_{M_{_{\alpha _{\eta }+1}}}-D_{M_{_{\alpha _{\eta
}}}}\right) \right\vert +\frac{M_{\alpha _{s}}^{1/p-1}}{\alpha _{s}}%
\left\vert \left( D_{j}-D_{M_{_{\alpha _{s}}}}\right) \right\vert  \notag \\
&\leq &\lambda \sum_{\eta =0}^{s-1}\frac{M_{\alpha _{\eta }}^{1/p}}{\alpha
	_{\eta }}+\frac{\lambda M_{\alpha _{s}}^{1/p}}{\alpha _{s}}\leq \frac{%
	2\lambda M_{\alpha _{s-1}}^{1/p}}{\alpha _{s-1}}+\frac{\lambda M_{\alpha
		_{s}}^{1/p}}{\alpha _{s}}\leq \frac{4\lambda M_{\alpha _{k-1}}^{1/p}}{\alpha
	_{k-1}}.  \notag
\end{eqnarray}

Let $M_{\alpha _{s-1}+1}+1\leq $ $j\leq M_{\alpha _{s}},$ where $s=1,...,k.$
Analogously to (\ref{8}) we can prove that

\begin{eqnarray*}
	&&\left\vert S_{j}f\right\vert =\left\vert \sum_{v=0}^{M_{\alpha _{s-1}+1}-1}%
	\widehat{f}(v)\psi _{v}\right\vert =\left\vert \sum_{\eta
		=0}^{s-1}\sum_{v=M_{\alpha _{\eta }}}^{M_{\alpha _{\eta }+1}-1}\frac{%
		M_{\alpha _{\eta }}^{1/p-1}}{\alpha _{\eta }}\psi _{v}\right\vert \\
	&=&\left\vert \sum_{\eta =0}^{s-1}\frac{M_{\alpha _{\eta }}^{1/p-1}}{\alpha
		_{\eta }}\left( D_{M_{_{\alpha _{\eta }+1}}}-D_{M_{_{\alpha _{\eta
	}}}}\right) \right\vert \leq \frac{2\lambda M_{\alpha _{s-1}}^{1/p}}{\alpha
		_{s-1}}\leq \frac{4\lambda M_{\alpha _{k-1}}^{1/p}}{\alpha _{k-1}}.
\end{eqnarray*}
Hence
\begin{equation}
\left\vert I\right\vert \leq \frac{1}{Q_{M_{\alpha _{k}}+2}}
\sum_{j=0}^{M_{\alpha _{k}}}q_{j}\left\vert S_{j}f\right\vert \leq \frac{4\lambda M_{\alpha _{k-1}}^{1/p}}{\alpha _{k-1}}\frac{1}{Q_{M_{\alpha
			_{k}}+2}}\sum_{j=0}^{M_{\alpha _{k}}}q_{j}\leq \frac{4\lambda M_{\alpha
		_{k-1}}^{1/p}}{\alpha _{k-1}}.  \label{10}
\end{equation}
If we now apply (\ref{6}) and (\ref{8}) we get that
\begin{eqnarray}
\left\vert II\right\vert &=&\frac{q_{M_{\alpha _{k}}+1}}{Q_{M_{\alpha_{k}}+2}}\left\vert \frac{%
	M_{\alpha _{k}}^{1/p-1}}{\alpha _{k}}\psi _{M_{\alpha _{k}}}+S_{M_{\alpha
		_{k}}}f\right\vert  \label{100} \\
&=&\frac{q_{M_{\alpha _{k}}+1}}{Q_{M_{\alpha_{k}}+2}}\left\vert \frac{M_{\alpha _{k}}^{1/p-1}}{%
	\alpha _{k}}\psi _{M_{\alpha _{k}}}+S_{M_{\alpha _{k-1}+1}}f\right\vert
\notag \\
&\geq &\frac{q_{M_{\alpha _{k}}+1}}{Q_{M_{\alpha_{k}}+2}}\left( \left\vert \frac{M_{\alpha
		_{k}}^{1/p-1}}{\alpha _{k}}\psi _{M_{\alpha _{k}}}\right\vert -\left\vert
S_{M_{\alpha _{k-1}+2}}f\right\vert \right)  \notag \\
&\geq &\frac{q_{M_{\alpha _{k}}+1}}{Q_{M_{\alpha _{k}}+2}}\left( \frac{M_{\alpha
		_{k}}^{1/p-1}}{\alpha _{k}}-\frac{4\lambda M_{\alpha _{k-1}}^{1/p}}{\alpha
	_{k-1}}\right) \notag\\ 
&\geq& \frac{q_{M_{\alpha _{k}}+1}}{Q_{M_{\alpha _{k}}+2}}\frac{M_{\alpha
		_{k}}^{1/p-1}}{4\alpha _{k}}.  \notag
\end{eqnarray}

Without lost the generality we may assume that $c=1$ in (\ref{cond1}). By
combining (\ref{10}) and (\ref{100}) we get
\begin{eqnarray}
	\left\vert T_{M_{\alpha _{k}}+2}f\right\vert &\geq &\left\vert II\right\vert
	-\left\vert I\right\vert \geq \frac{q_{M_{\alpha _{k}}+1}}{Q_{M_{\alpha _{k}}+2}}\frac{	M_{\alpha _{k}}^{1/p-1}}{4\alpha _{k}}-\frac{4\lambda M_{\alpha _{k-1}}^{1/p}}{\alpha _{k-1}}\label{155aba2} \\
	&\geq &\frac{M_{\alpha _{k}}^{1/p-2}}{4\alpha _{k}}-\frac{4\lambda M_{\alpha
			_{k-1}}^{1/p}}{\alpha _{k-1}}\geq \frac{M_{\alpha _{k}}^{1/p-2}}{16\alpha _{k}}. \notag
\end{eqnarray}
On the other hand,
\begin{equation}
\mu \left\{ x\in G_{m}:\left\vert T_{M_{\alpha _{k}}+2}f\left( x\right)
\right\vert \geq \frac{M_{\alpha _{k}}^{1/p-2}}{16\alpha _{k}}\right\} =\mu
\left( G_{m}\right) =1.  \label{88}
\end{equation}
Let $0<p<1/2.$ Then
\begin{eqnarray}
&&\frac{M_{\alpha _{k}}^{1/p-2}}{16\alpha _{k}}\cdot \left(\mu \left\{ x\in
G_{m}:\left\vert T_{M_{\alpha _{k}}+2}f\left( x\right) \right\vert \geq
\frac{M_{\alpha _{k}}^{1/p-2}}{16\alpha _{k}}\right\} \right)^{1/p} \label{99} \\
&=&\frac{M_{\alpha _{k}}^{1/p-2}}{16\alpha _{k}}\rightarrow \infty ,\text{ \	as }k\rightarrow \infty .  \notag
\end{eqnarray}%
\textbf{\ }
The proof is complete.

\textbf{Proof of Corollary \ref{cor0}.} 
Since $R_{n},	U_{n}^{\alpha} \text{ and } V_{n}^{\alpha}$ are  the $T$ means with non-increasing sequence $\{q_{k}:k\geq 0\},$ then the proof of this corollary is direct consequence of Theorem \ref{theorem1}.

\textbf{Proof of Corollary \ref{cor2}.} 
According to Theorem 1 a) and Proposition \ref{prop2} we also have weak $(1,1)$ type inequality and by well-known density argument due to Marcinkiewicz and Zygmund \cite{mz} we have $T_{n}f\rightarrow f,$ a.e., for all $f\in L_1.$ Which follows proof of Corollary \ref{cor2}.

\textbf{Proof of Corollary \ref{cor1}.} 
Since $R_{n},	U_{n}^{\alpha} \text{ and } V_{n}^{\alpha}$ are  the $T$ means with non-increasing sequence $\{q_{k}:k\geq 0\},$ then the proof of this corollary is direct consequence of Corollary \ref{cor2}.

\textbf{Proof of Theorem \ref{theorem2}.} 
Let the sequence $\{q_{k}:k\geq 0\}$ be
non-decreasing. By combining (\ref{2b}) with (\ref{2c}) and using Abel transformation we get that
\begin{eqnarray*}
	\left\vert T_{n}f\right\vert &\leq &\left\vert \frac{1}{Q_{n}}\overset{n-1}{%
		\underset{j=1}{\sum }}q_{j}S_{j}f\right\vert \\
	&\leq &\frac{1}{Q_{n}}\left( \overset{n-2}{\underset{j=1}{\sum }}\left\vert
	q_{j}-q_{j+1}\right\vert j\left\vert \sigma _{j}f\right\vert
	+q_{n-1}(n-1)\left\vert \sigma _{n}f\right\vert \right) \\
	&\leq &\frac{1}{Q_{n}}\left( \overset{n-2}{\underset{j=1}{\sum }}-\left(
	q_{j}-q_{j+1}\right) j-q_{n-1}(n-1)+2q_{n-1}(n-1)\right) \sigma ^{\ast }f \\
	&\leq &\frac{1}{Q_{n}}\left( 2q_{n-1}(n-1)-Q_{n}\right) \sigma ^{\ast }f\leq c\sigma ^{\ast	}f
\end{eqnarray*}
so that
\begin{equation}
	T^{\ast }f\leq c \sigma ^{\ast }f.  \label{12aaaaa}
\end{equation}

If we apply (\ref{12aaaaa}) and Theorem W1 we can conclude that the maximal
operators $T^{\ast }$ are bounded from the Hardy space $H_{1/2}$ to the
space $weak-L_{1/2}.$ 
The proof of part a) is complete.

To prove part b) of Theorem 2 we use the martingale, defined by (\ref{55}) where $\alpha_k$ satisfy conditions (\ref{2})-(\ref{4}). It is easy to show that for every non-increasing sequence $\{q_{k}:k\geq 0\}$ it automatically holds that
\begin{equation*}
q_{M_{\alpha _{k}+1}}/Q_{M_{\alpha _{k}+2}}\geq c/M_{\alpha _{k}}.
\end{equation*}

According to (\ref{155aba})-(\ref{155aba2}) we can conclude that
\begin{equation*}
\left\vert T_{M_{\alpha _{k}}+2}f\right\vert \geq \left\vert II\right\vert
-\left\vert I\right\vert \geq \frac{M_{\alpha _{k}}^{1/p-2}}{8\alpha _{k}}.
\end{equation*}

Analogously to\ (\ref{88}) and (\ref{99}) we then get that%
\begin{equation*}
\sup_{k\in \mathbb{N}}\left\Vert T_{M_{\alpha _{k}}+2}f\right\Vert _{weak-L_{p}}=\infty .
\end{equation*}

The proof is complete.

\textbf{Proof of Corollary \ref{cor4}.} Since $B^{\alpha,\beta, *}$ are  the $T$ means with non-decreasing sequence $\{q_{k}:k\geq 0\},$ then the proof of this corollary is direct consequence of Theorem \ref{theorem2}.

 \textbf{Proof of Corollary \ref{cor5}.} According to Proposition \ref{prop2} we can conclude that $T^{\ast }$ has weak type-(1,1) and by well-known density argument due to Marcinkiewicz and Zygmund \cite{mz} we also have $T_{n}f\rightarrow f,$ a.e.. Which follows proof of Corollary \ref{cor5}.

\textbf{Proof of Corollary \ref{cor6}.} Since $B^{\alpha,\beta, *}$ are  the $T$ means with non-decreasing sequence $\{q_{k}:k\geq 0\},$ then the proof of this corollary is direct consequence of Corollary  \ref{cor5}.

\end{document}